\documentclass{article}[12pt]
\usepackage{amssymb , amsmath}

\newtheorem{theorem}{\sc Theorem}[section]
\newtheorem{lemma}{\sc Lemma}[section]

\newcommand{\be}{ \begin{equation} }
\newcommand{\ee}{ \end{equation} }

\newcommand{\hh}{\widehat{h}}
\newcommand{\hu}{\widehat{u}}

\def\C    {\hbox{\kern 2.7pt\vrule height6pt depth-.35pt \rm \kern-2.7pt C}}

\title{On the resolvent technique for stability of plane Couette
flow}
\author{Pablo Braz e Silva\footnotemark[2]\ \footnotemark[3]}

\begin{document}
\date{}
\maketitle
\renewcommand{\thefootnote}{\fnsymbol{footnote}}

\footnotetext[3]{Supported by a post-doctoral
       fellowship FAPESP/Brazil: 02/13270-1}
\footnotetext[2]{Instituto de Matem\'atica, Estat\'istica e
Computa\c c\~ao Cient\'ifica - UNICAMP, Cx. Postal 6065, CEP
13083-970, Campinas, SP, Brazil ({\tt pablo@ime.unicamp.br}).}

\renewcommand{\thefootnote}{\arabic{footnote}}

\begin{abstract}
 We discuss the application of the resolvent technique to prove stability of
 plane Couette flow.
 Using this technique, we derive a threshold amplitude for perturbations that
 can lead to turbulence in terms of the Reynolds number. Our main
 objective is to show exactly how much control
 one should have over the perturbation to assure stability via this technique.
\end{abstract}

\thanks{ {\em 2000 Mathematics Subject Classifications:} 76E05, 47A10, 35Q30, 76D05.
\hfil\break\noindent {\em Key words and phrases:} Couette flow,
resolvent estimates

\section{Introduction}

   We discuss stability of plane Couette flow via the resolvent
   method. Applying this method, one can derive lower bounds for
   the norms of perturbations that can lead to turbulence.
   Our aim is to discuss and clarify a point that has been overseen so far,
   which is
   to determine how much control over the perturbations
   one should assume to derive the stability result via the resolvent technique.
   We discuss the two dimensional case, but all the considerations can be used for the
   three dimensional case with minor technical changes.
   The main difference between the two and three spatial dimensions
   is that different resolvent estimates
   hold for each case, leading to different thresholds. This point will be made clear later
   on.
   We begin describing the problem and discussing previous works using
   the resolvent method.

\section{The problem and known results}

   We are interested in the following initial boundary value problem:
\begin{equation}   \label{eq1}
    \left\{ \begin{array}{l}
        \displaystyle u_t + (u\cdot \nabla )u + \nabla p = \frac{1}{R}\Delta u \\
        \displaystyle \nabla \cdot u = 0 \\
        u(x,0,t) = (0,0) \\
        u(x,1,t) = (1,0) \\
        u(x,y,t) = u(x+1,y,t) \\
        u(x,y,0) = f(x,y)
     \end{array} \right.
\end{equation}
where $u : \mathbb{R} \times [0 , 1] \times [0 , \infty)
\longrightarrow \mathbb{R}^2$ is the unknown function $u(x,y,t) =
(u_1(x,y,t),u_2(x,y,t))$. The positive parameter $R$ is the
Reynolds number. The initial condition $f(x,y)$ is assumed to be
smooth, divergence free and compatible with the boundary
conditions. The pressure $p(x,y,t)$ can be determined in terms of
$u$ by the elliptic problem
\begin{equation}   \label{eq2}
    \left \{ \begin{array}{l}
        \displaystyle \Delta p = -\nabla \cdot ((u\cdot \nabla )u) \\
        \displaystyle p_y(x,0,t) = \frac{1}{R} u_{2yy}(x,0,t) \vspace{.1cm}\\
        \displaystyle p_y(x,1,t) = \frac{1}{R} u_{2yy}(x,1,t) .
     \end{array} \right.
\end{equation}
It can be easily seen that $U(x,y) = (y,0)$, $P =  constant$ is a
steady solution of problem (\ref{eq1}). The vector field $U(x,y) =
(y,0)$ is known as Couette flow.

Using the resolvent technique, one can prove and quantify
asymptotic stability for this flow. By quantification we mean the
derivation a number $M(R)$ such that disturbances of the flow with
norm less than $M(R)$ will tend to zero as time $t$ tends to
infinity. In other words, deriving a lower bound for the norm of
perturbations that can lead to turbulence. For general discussion
about the resolvent technique, see \cite{KL2} and \cite{KL3}.

This problem has been studied for the 3 spatial dimensions case in
\cite{K}, and a threshold amplitude for perturbations was found to
be of order $\mathcal{O}
$$(R^{-\frac{21}{4}})$. The estimates of the resolvent of the
linearized equations governing perturbations where those found
numerically in \cite{RE} and \cite{T}, predicting the resolvent
constant of the linear operator associated with the problem to be
proportional to $R^2$. In \cite{L1}, the resolvent technique was
used again to prove the stability of the 3 dimensional problem but
the estimates for the resolvent constant were those in \cite{L2}.
By using modified norms, the authors achieve $M(R)$ of order
$\mathcal{O}$$(R^{-3})$ for two of the components of the
perturbation, and of order $\mathcal{O}$$(R^{-4})$ for the
remaining component. Our approach uses again the resolvent
technique, and we use the same norms as \cite{K}, with the obvious
modifications for the 2 spatial dimensions case. We show that this
approach leads, in our case, to a threshold amplitude of order
$\mathcal{O}$$(R^{-3})$. We note that our argument is the same
used in \cite{K}, with some minor differences. The only reason for
the better exponent in our case is the better dependence of the
resolvent constant on $R$ for the 2 dimensional case. In this
case, the resolvent constant is proportional to $R$, as found in
\cite{BS}. We carry out the argument in details again only because
it is important for our aim, which is to clarify a subtle point
that has been overseen in previous works: In \cite{K}, it was said
that one needs
 control over the sobolev norm $H^2$ of the perturbation to assure stability. Later
 on, in \cite{L1}, the authors note that the
 $H^2$ is not enough, and claim that one needs control over the
 norm $H^4$. Actually, this is not enough yet, since in one of the
 directions, one needs control over six derivatives of the perturbation.
 This necessity is due to the pressure terms appearing in the problem.
 In section \ref{section4}, we show
 in details estimates for these terms, and clarify
 the reason of this requirement. Moreover, our argument shows that
 derivatives of different orders of the perturbation scale differently
 with the Reynolds number. In other words, to assure decay of the perturbations
 via the resolvent method, one should require the perturbation to be small
 in some weighted norm involving six derivatives, where the weights depend
 on the Reynolds number $R$.

This work is divided in 4 sections: in section \ref{section2}, we
introduce some basic notation and derive the equations for
perturbations of the Couette flow; in section \ref{section3} we
derive estimates for the solution of the linearized equations for
the perturbations; in \ref{section4}, we use those estimates to
prove asymptotic stability for the flow, and to derive the
threshold amplitude $M(R)$. In section \ref{section5}, we derive
carefully the estimates for the pressure terms involved in the
problem.

\section{Notation and equations for the
perturbations}\label{section2}

 We denote by $\langle \cdot , \cdot \rangle$ and
$\|\cdot\|$ the $L_2$ inner product and norm over $\Omega = [0 ,1]
\times [0 ,1]$:
$$\begin{array}{lcr}
      \displaystyle \langle u , w \rangle = \int_{\Omega} \overline{u} \cdot w \hspace{.1cm} dxdy
      & ; &
      \displaystyle \|u\|^2 = \langle u , u \rangle .
   \end{array}
$$
All the matrix norms that appear in this paper are the usual
Frobenius norms. The usual sobolev norm of $u$ over $\Omega$ is
denoted by
$$ \displaystyle \|u\|^2_{H^n (\Omega )} = \sum_{j=0}^{n}\|D^j u\|^2$$
where $D^j$ denotes the $j$-th derivative of $u$ with respect to
the space variables. Unless stated otherwise, all norms in the
space variables will be calculated over $\Omega$ and therefore we
will write $\displaystyle \| \cdot \|_{H^n(\Omega)}$ as
$\displaystyle \| \cdot \|_{H^n}$. We make use of a 2 dimensional
version of the weighted norm $\|\cdot\|_{\widetilde{H}}$ used in
Kreiss\cite{K}:
\begin{equation}\label{eq3}
 \displaystyle \|u\|^2_{\widetilde{H}} = \|u\|^2 + \frac{1}{R}\|D u\|^2 +
\frac{1}{R^2} \|u_{xy}\|^2 .
\end{equation}
We also define another weighted norm $\|\cdot\|_{H^6_m}$ by
\begin{eqnarray}
 \displaystyle \|u\|^2_{H^6_m} & = & \|u\|^2_{H^2} + \frac{1}{R^2}\|D^3 u\|^2 + \frac{1}{R^2}\|D^4 u\|^2 \nonumber \\
 & & \mbox{} + \frac{1}{R^4}\| u_{2xxyyy}\|^2 +\frac{1}{R^4}\| u_{2yyyyy}\|^2 +
\frac{1}{R^4}\| u_{2yyyyyy}\|^2,
\end{eqnarray}
where $u = (u_1 , u_2 )$.

The maximum norm over $\Omega$ is denoted by $\displaystyle |\cdot
|_{\infty}$. The norm $\displaystyle \|\cdot\|_{\widetilde{H}}$ is
related with the maximum norm by the sobolev type inequality (see
\cite{KL}, Appendix 3, Theorem A.3.14)
$$
\displaystyle |\cdot|_{\infty}^2 \leq \widetilde{C} R
\|\cdot\|_{\widetilde{H}}^2 .
$$
Since we are interested in functions which are also dependent on
time, we use that
\begin{equation} \label{eq4}
\displaystyle |u(\cdot , t) |_{\infty}^2 \leq \widetilde{C} R \|u(
\cdot, t )\|_{\widetilde{H}}^2 \; , \; \forall \: t \geq 0,
\end{equation}
where $\widetilde{C}$ is a constant independent of any of the
parameters.

We are interested first in proving asymptotic stability for the
Couette flow, which is a  stationary solution of (\ref{eq1}), that
is, to prove that perturbations of the stationary solution that
are small enough in some norm will tend to $0$ as $t$ tends to
infinity. More specifically, we will show that perturbations
having norm $\| \cdot \|_{H^6_m}$ of order $R^{-3}$ decay with
time.

To this end, let $U = U(x,y)$, $P = P(x,y)$ be a stationary
solution of (\ref{eq1}). We can obviously use the Couette flow,
but we think that the structure of the argument is easier to be
understood if one uses any stationary solution. This will not
change the estimates we will prove. We derive the equations
satisfied by perturbations of this base flow. Let $u(x,y,t)$ ,
$p(x,y,t)$ be a solution of (\ref{eq1}) with initial condition
$f(x,y) = U(x,y) + \epsilon f'(x,y)$, where $f^\prime$ is
divergence free and $\|f^\prime \|_{H^6_m(\Omega)} =1$. Then,
$\epsilon$ defines a unique perturbation amplitude. Write
$u(x,y,t) = U(x,y) + \epsilon u'(x,y,t)$ and $p(x,y,t) = P(x,y) +
\epsilon p'_1(x,y,t) +\epsilon^2 p'_2(x,y,t)$. Then
$u^\prime$,$p'_1$,$p'_2$ satisfy the system
$$
    \left\{ \begin{array}{l}
        \displaystyle u'_t + (u'\cdot \nabla )U + (U\cdot \nabla )u'+ \nabla p'_1 +
\epsilon (u'\cdot \nabla )u' + \epsilon \nabla p'_2 = \frac{1}{R}\Delta u' \\
        \displaystyle \nabla \cdot u' = 0 \\
        u'(x,0,t) = (0,0) \\
        u'(x,1,t) = (0,0) \\
        u'(x,y,t) = u'(x+1,y,t) \\
        u'(x,y,0) = f'(x,y).
     \end{array} \right.
$$
The functions $p'_1$ and $p'_2$ are given in terms of $u'$ by
$$
    \left \{ \begin{array}{l}
        \displaystyle \Delta p'_1 = -\nabla \cdot ((u'\cdot \nabla )U) -  \nabla \cdot ((U\cdot \nabla )u') \\
        \displaystyle p'_{1y}(x,0,t) = \frac{1}{R} u'_{2yy}(x,0,t) \vspace{.1cm}\\
        \displaystyle p'_{1y}(x,1,t) = \frac{1}{R} u'_{2yy}(x,1,t)
     \end{array} \right.
$$
and
$$
    \left \{ \begin{array}{l}
        \displaystyle \Delta p'_2 = -\nabla \cdot ((u'\cdot \nabla )u') \\
        \displaystyle p'_{2y}(x,0,t) = 0\vspace{.1cm}\\
        \displaystyle p'_{2y}(x,1,t) = 0 .
     \end{array} \right.
$$
As we show in section \ref{section5}, the functions $p'_1$ and
$p'_2$ can be estimated in terms of $u'$ by
\begin{eqnarray*}
\|\nabla p'_1 (\cdot , \cdot , t )\|^2 & \leq & C \left( \|u'
(\cdot , \cdot , t )\|_{H^1}^2 + \frac{1}{R^2} \| u'_{2yy}(\cdot ,
\cdot , t ) \| + \frac{1}{R^2} \| u'_{2yyy} (\cdot , \cdot , t )\|
\right)
   \, , \, \forall \, t\geq 0 ,\\
\|\nabla p'_2 (\cdot , \cdot , t )\|^2 & \leq & \|( u'\cdot
\nabla) u'(\cdot , \cdot , t )\|^2
   \, , \, \forall \, t\geq 0.
\end{eqnarray*}
From now on, to simplify the notation, we drop the $'$ in the
equations above, and just write $u$, $p_1$, $p_2$. With this
notation, the equations above are
\begin{equation}   \label{eq5}
    \left\{ \begin{array}{l}
        \displaystyle u_t + (u\cdot \nabla )U + (U\cdot \nabla )u+ \nabla p_1 +
\epsilon (u\cdot \nabla )u + \epsilon \nabla p_2 = \frac{1}{R}\Delta u \\
        \displaystyle \nabla \cdot u = 0 \\
        u(x,0,t) = (0,0) \\
        u(x,1,t) = (0,0) \\
        u(x,y,t) = u(x+1,y,t) \\
        u(x,y,0) = f(x,y),
     \end{array} \right.
\end{equation}
\begin{equation}   \label{eq6}
    \left \{ \begin{array}{l}
        \displaystyle \Delta p_1 = -\nabla \cdot ((u\cdot \nabla )U) -  \nabla \cdot ((U\cdot \nabla )u) \\
        \displaystyle p_{1y}(x,0,t) = \frac{1}{R} u_{2yy}(x,0,t) \vspace{.1cm}\\
        \displaystyle p_{1y}(x,1,t) = \frac{1}{R} u_{2yy}(x,1,t)
     \end{array} \right.
\end{equation}
and
\begin{equation}   \label{eq7}
    \left \{ \begin{array}{l}
        \displaystyle \Delta p_2 = -\nabla \cdot ((u\cdot \nabla )u) \\
        \displaystyle p_{2y}(x,0,t) = 0\vspace{.1cm}\\
        \displaystyle p_{2y}(x,1,t) = 0 .
     \end{array} \right.
\end{equation}
Note that $p_1$ depends linearly on $u$. Moreover, for all $t\geq
0$, we have
\begin{eqnarray}
\|\nabla p_1 (\cdot , \cdot , t )\|^2 & \leq & C \left( \|u (\cdot
, \cdot , t )\|_{H^1}^2 + \frac{1}{R^2} \| u_{2yy}(\cdot , \cdot ,
t ) \| + \frac{1}{R^2} \| u_{2yyy} (\cdot , \cdot , t )\| \right),
    \label{eqp1}\\
\|\nabla p_2 (\cdot , \cdot , t )\|^2 & \leq & \|( u\cdot \nabla)
u(\cdot , \cdot , t )\|^2 . \label{eqp2}
\end{eqnarray}
When the initial data is divergence free and the terms of pressure
are given by the equations (\ref{eq6}) and (\ref{eq7}) above, the
solution $u$ of problem (\ref{eq5}) remains divergence free for
all time $t$. Therefore, we drop the continuity equation and write
problem (\ref{eq5}) as
\begin{equation}   \label{eq8}
    \left\{ \begin{array}{l}
        \displaystyle u_t  = {\mathcal L} u - \epsilon (u\cdot \nabla )u - \epsilon \nabla p_2 \\
        u(x,0,t) = (0,0) \\
        u(x,1,t) = (0,0) \\
        u(x,y,t) = u(x+1,y,t) \\
        u(x,y,0) = f(x,y),
     \end{array} \right.
\end{equation}
where ${\mathcal L}$ is a linear operator depending on the
parameter $R$, defined by
\begin{equation}\label{operator}
\displaystyle {\mathcal L} u = \frac{1}{R}\Delta u - (u\cdot
\nabla )U - (U\cdot \nabla )u - \nabla p_1 ,
\end{equation}
with $p_1$ given by (\ref{eq6}). It is very important to note that
this linear operator has also an integral part, which is the term
$\nabla p_1$. Moreover, as inequality (\ref{eqp1}) shows, to
estimate $\| \nabla p_1 \|$ one needs three space derivatives of
the second component of $u$, at least in one of the directions.

This integral part of the operator $\mathcal L$ seems to be the
point which was so far overseen. If one neglects that and consider
only to the differential part of $\mathcal L$, one will be led to
conclude that $\| \cdot \|_{H^2}$ will be enough to estimate the
right hand side of (\ref{operator}).

We first apply the resolvent technique to prove stability of the
stationary flow. For that end, it is convenient to have
homogeneous initial conditions. Therefore, we transform the
problem (\ref{eq8}) to a similar problem with homogeneous initial
condition by defining
\begin{equation} \label{eq9}
   v(x,y,t) := u(x,y,t) - e^{-t}f(x,y).
\end{equation}
Note that $v$ and $u$ have the same behavior as $t \rightarrow
\infty$. Moreover, $v$ given by (\ref{eq9}) satisfies
\begin{equation}   \label{eq10}
    \left\{ \begin{array}{l}
        \displaystyle v_t  = {\mathcal L} v - \epsilon \{ (v\cdot \nabla )v + e^{-t}(v\cdot \nabla) f +
e^{-t}(f\cdot \nabla) v \} - \epsilon \nabla p_2 + F(x,y,t)   \\
        v(x,0,t) = (0,0) \\
        v(x,1,t) = (0,0) \\
        v(x,y,t) = v(x+1,y,t) \\
        v(x,y,0) = (0,0),
     \end{array} \right.
\end{equation}
where $F(x,y,t) = e^{-t}(({\mathcal L} + I )f -\epsilon e^{-t} (f
\cdot \nabla)f)$. Note that $F$, $F_t \in
L_2([0,\infty);L_2(\Omega))$, that is, both $\|F(\cdot , \cdot
,t)\|^2$ and $\|F_t(\cdot , \cdot ,t)\|^2$ are integrable over
$[0,\infty)$.

\section{Linear Problem}\label{section3}

We first consider the general linear problem
\begin{equation}   \label{eq11}
    \left\{ \begin{array}{l}
        \displaystyle v_t  = {\mathcal L} v + F(x,y,t)   \\
        v(x,0,t) = (0,0) \\
        v(x,1,t) = (0,0) \\
        v(x,y,t) = v(x+1,y,t) \\
        v(x,y,0) = (0,0),
     \end{array} \right.
\end{equation}
where $\|F(\cdot , t) \|^2$ and $\|F_t (\cdot , t) \|^2$
integrable over the domain $[0,\infty)$:
$$
   \displaystyle \int_0^{\infty} \left(\|F(\cdot , t) \|^2 + \|F_t (\cdot , t) \|^2
\right) dt < \infty.
$$
In our case of two spatial dimensions, resolvent estimates were
found in \cite{BS}:
\begin{equation} \label{eq12}
  \| \widetilde{v}(\cdot , s) \|^2 \leq C_1 R^2
\| \widetilde{F} (\cdot , s)\|^2 \hspace{.2cm}, \hspace{.2cm}
\mbox{Re} s \geq 0 ,
\end{equation}
where \, $\widetilde{}$ \, stands for the Laplace transform with
respect to $t$ , $s$ is its variable and $C_1$ is an absolute
constant, that is, it does not depend on any of the parameters or
functions. One can prove, as in \cite{K}, Appendix A, that
(\ref{eq12}) implies
\begin{equation} \label{eq13}
  \| \widetilde{v}(\cdot , s) \|^2_{\widetilde{H}} \leq C R^2 \| \tilde{F}(\cdot , s) \|^2
\end{equation}
where $C$ depends on $C_1$ and on $U$ and its first derivative.
Since for our problem $U$ is fixed as the Couette flow, $C$ is an
absolute constant as well. From now on, we will use $C$ for any
absolute constant, and replace its value as necessary keeping the
notation $C$. No attempt is made to optimize those constants,
since the most important result is the dependence of the threshold
amplitude on the Reynolds number.

Using Parseval's relation, inequality (\ref{eq13}) for the
transformed functions is translated to the original functions as
\begin{equation} \label{eq14}
  \int_0^{\infty}  \| v(\cdot ,t) \|^2_{\widetilde{H}}dt \leq
   C R^2 \int_0^{\infty} \| F(\cdot , t) \|^2 dt.
\end{equation}
Obviously, $\displaystyle \int_0^{T}  \| v(\cdot ,t)
\|^2_{\widetilde{H}}dt \leq \int_0^{\infty}  \| v(\cdot ,t)
\|^2_{\widetilde{H}}dt $ , $\displaystyle \forall \, T \geq 0$.
Moreover, since the solution of the equation up to time $T$ does
not depend on the forcing $F(x,y,t)$ for $t>T$, we have
\begin{equation} \label{eq15}
   \begin{array}{lcr}
      \displaystyle \int_0^{T}  \| v(\cdot ,t) \|^2_{\widetilde{H}}dt \leq
      C R^2 \int_0^{T} \| F(\cdot , t) \|^2 dt & ,  & \forall \, T \geq 0.
    \end{array}
\end{equation}
For our argument, we also need similar estimates for $v_t$. To
this end, differentiate equation (\ref{eq11}) to get
\begin{equation}   \label{eq16}
    \left\{ \begin{array}{l}
        \displaystyle v_{tt}  = {\mathcal L} v_t + F_t(x,y,t)   \\
        v_t(x,0,t) = (0,0) \\
        v_t(x,1,t) = (0,0) \\
        v_t(x,y,t) = v_t(x+1,y,t) \\
        v_t(x,y,0) = F(x,y,0) =: g(x,y),
     \end{array} \right.
\end{equation}
that is, $v_t$ satisfies an equation of the same type as
(\ref{eq11}), but with non-homogeneus initial conditions $g(x,y) =
F(x,y,0)$. Performing the same type of initialization as before,
that is, defining $\varphi := v_t - e^{-t} g $, we get a similar
problem for $\varphi$, with homogeneus initial conditions and an
extra forcing term. Using the estimates for the resolvent, and
writing those in terms of $v_t$, we get
\begin{equation} \label{eq17}
   \begin{array}{lll}
      \displaystyle \int_0^{T}  \| v_t (\cdot ,t) \|^2_{\widetilde{H}}dt  & \leq &
      \displaystyle \| F(x,y,0) \|^2_{\widetilde{H}} + C R^2 \| ({\mathcal L} +I)F(x,y,0) \|^2  \\
       & &  \displaystyle\mbox{} + C R^2 \int_0^{T} \| F_t (\cdot , t) \|^2 dt \hspace{.2cm} , \hspace{.2cm}
          \forall T \geq 0.
    \end{array}
\end{equation}
Combining (\ref{eq15}) and (\ref{eq17}) gives, for $v$ the
solution of (\ref{eq11}),
\begin{equation} \label{eq18}
   \begin{array}{l}
       \displaystyle \int_0^{T}
\left( \| v(\cdot ,t) \|^2_{\widetilde{H}} + \| v_t(\cdot ,t)
\|^2_{\widetilde{H}}\right) dt
       \leq
        \| F(x,y,0) \|^2_{\widetilde{H}} + C R^2 \| ({\mathcal L} +I)F(x,y,0) \|^2  \\
         \displaystyle \mbox{} +
          C R^2 \int_0^{T} \left( \| F (\cdot , t) \|^2 + \| F_t (\cdot , t) \|^2\right) dt
        \hspace{.3cm} , \hspace{.3cm}\forall T \geq 0.
    \end{array}
\end{equation}
Now, using these estimates for the solution of the linear problem,
we can prove a stability result for the nonlinear equation.

\section{Stability for the Nonlinear Problem}\label{section4}

The nonlinear problem (\ref{eq10}) is
\begin{equation}   \label{eq19}
    \left\{ \begin{array}{l}
        \displaystyle v_t  = {\mathcal L} v - \epsilon \{ (v\cdot \nabla )v + e^{-t}(v\cdot \nabla) f +
e^{-t}(f\cdot \nabla) v \} - \epsilon \nabla p_2 + F(x,y,t)   \\
        v(x,0,t) = (0,0) \\
        v(x,1,t) = (0,0) \\
        v(x,y,t) = v(x+1,y,t) \\
        v(x,y,0) = (0,0),
     \end{array} \right.
\end{equation}
where $F(x,y,t) = e^{-t}(({\mathcal L} + I )f -\epsilon e^{-t} (f
\cdot \nabla)f)$. We prove the following:
\begin{theorem}
There exists $\epsilon_0 > 0$, $\epsilon_0 = \epsilon_0 (R)$, such
that if $0\leq |\epsilon | < \epsilon_0$, then the solution
$v(x,y,t)$ of (\ref{eq19}) satisfies
$$ \displaystyle
\lim_{t\rightarrow \infty} |v(\cdot , t)|_{\infty} = 0 .
$$
Moreover, $\epsilon_0 = {\mathcal O}(R^{-3})$.
\end{theorem}
\paragraph{Proof:}
We consider problem (\ref{eq19}) as a linear problem with forcing
\begin{equation} \label{eq20}
G(x,y,t) := F(x,y,t) - \epsilon \{ (v\cdot \nabla )v +
e^{-t}(v\cdot \nabla) f + e^{-t}(f\cdot \nabla) v \} - \epsilon
\nabla p_2.
\end{equation}
Applying inequality (\ref{eq18}) with forcing term $G$ gives
\begin{eqnarray} \label{eq21}
    \displaystyle \int_0^{T}  \left( \| v(\cdot ,t) \|^2_{\widetilde{H}} +
\| v_t(\cdot ,t) \|^2_{\widetilde{H}}\right) dt
       \leq
       \displaystyle \| G(x,y,0) \|^2_{\widetilde{H}} + C R^2 \| ({\mathcal L} +I)G(x,y,0) \|^2 +  \nonumber \\
        \displaystyle \mbox{} + C R^2 \int_0^{T} \left( \| G (\cdot , t) \|^2 + \| G_t (\cdot , t) \|^2\right) dt
        \hspace{.3cm} \forall T \geq 0.
    \end{eqnarray}
From the definition of $G$, we have
\begin{equation}\label{equacao1}
\begin{array}{l}
   \displaystyle \int_0^{T}  \left( \| v(\cdot ,t) \|^2_{\widetilde{H}} +
           \| v_t(\cdot ,t) \|^2_{\widetilde{H}}\right) dt   \leq 2\| F(x,y,0) \|^2_{\widetilde{H}}
           + 2\epsilon^2 \|\nabla p_2(x,y,0)\|^2_{\widetilde{H}} \\
   + C R^2 \| ({\mathcal L}_R +\mathcal{I})F(x,y,0) \|^2 + CR^2 \|({\mathcal L}_R +\mathcal{I})p_2(x,y,0) \|^2\\
    \displaystyle \mbox{} +  C R^2 \int_0^{T} \left( \| F - \epsilon \{ (v\cdot \nabla )v +
   e^{-t}(v\cdot \nabla) f + e^{-t}(f\cdot \nabla) v \} - \epsilon \nabla p_2 \|^2\right) dt  \\
    \displaystyle \mbox{} +  C R^2\int_0^T \left( \| (F - \epsilon \{ (v\cdot \nabla )v + e^{-t}(v\cdot \nabla) f +
e^{-t}(f\cdot \nabla) v \} - \epsilon \nabla p_2)_t \|^2 \right)
dt.
\end{array}
\end{equation}
Since $p_2$ is given by (\ref{eq7}), we have (see section
\ref{section5})
$$
\begin{array}{lcr}
  \displaystyle \|\nabla p_2 \| \leq \|( u\cdot \nabla) u\| & ; &
   \displaystyle \|(\nabla p_2)_t \| \leq \|(( u\cdot \nabla) u)_t\|.
\end{array}
$$ Thus, using (\ref{eq9}), we can estimate $\nabla p_2$ by $f$ and $v$.
Moreover,
$$\displaystyle \|\nabla p_2(\cdot,\cdot,0) \|^2
\leq \|( u\cdot \nabla) u(\cdot,\cdot,0)\|^2 = \|( f\cdot \nabla)
f\|^2,$$ and since $\| f \|^2_{H^6_m} = 1$, inequality
(\ref{equacao1}) gives
$$
\begin{array}{l}
\displaystyle \int_0^{T}  \left( \| v(\cdot ,t)
\|^2_{\widetilde{H}} +
           \| v_t(\cdot ,t) \|^2_{\widetilde{H}}\right) dt   \leq
 \| F(x,y,0) \|^2_{\widetilde{H}} + C R^2 \| ({\mathcal L} +I)F(x,y,0) \|^2 \\
\displaystyle \mbox{} + C R^2 \int_0^{\infty} \left( \| F \|^2 +
\| F_t \|^2 \right) dt \displaystyle \mbox{} + C R^2 \epsilon^2
\int_0^{T} \left( \| (v\cdot \nabla )v \|^2 +
 \| (v_t \cdot \nabla )v \|^2 + \| (v\cdot \nabla )v_t \|^2 \right) dt \\
 \displaystyle \mbox{} + C R^2 \epsilon^2\int_0^T \left( \| e^{-t}(v\cdot \nabla) f\|^2 + \| e^{-t}(f\cdot \nabla) v \|^2
 + \| e^{-t}(v_t \cdot \nabla) f \|^2 + \| e^{-t}(f\cdot \nabla) v_t \|^2\right)  dt.
\end{array}
$$
Since $$F(x,y,t) = e^{-t}(({\mathcal L} + I )f -\epsilon e^{-t} (f
\cdot \nabla)f),$$ we have $F(x,y,0) = ({\mathcal L} + I )f
-\epsilon (f \cdot \nabla)f:={\mathcal P}f$. With this notation,
the inequality above is
\begin{equation}\label{ineq1}
\begin{array}{l}
  \displaystyle \int_0^{T}  \left( \| v(\cdot ,t) \|^2_{\widetilde{H}} +
\| v_t(\cdot ,t) \|^2_{\widetilde{H}}\right) dt  \leq
         \| {\mathcal P} f \|^2_{\widetilde{H}} + C R^2 \| ({\mathcal L} +I){\mathcal P} f \|^2
+C R^2 \| ({\mathcal L} +I)f\|^2 \\
\displaystyle \mbox{} + C R^2 \epsilon^2 \|(f\cdot \nabla)f\|^2 +
C R^2 \epsilon^2 \int_0^{T} \left( \| (v\cdot \nabla )v \|^2 +
 \| (v_t \cdot \nabla )v \|^2 + \| (v\cdot \nabla )v_t \|^2\right) dt \\
 \displaystyle \mbox{}
+ C R^2 \epsilon^2 \int_0^T \left( \| e^{-t}(v\cdot \nabla) f\|^2
+ \| e^{-t}(f\cdot \nabla) v \|^2  + \| e^{-t}(v_t \cdot \nabla) f
\|^2 +
        + \| e^{-t}(f\cdot \nabla) v_t \|^2\right)  dt .
\end{array}
\end{equation}
To apply the resolvent method, one needs control over the terms
depending on $f$ of the right hand side of inequality
(\ref{ineq1}). Its second term is
$$
\| ({\mathcal L} +I){\mathcal P} f\| = \|({\mathcal L}
+I)(({\mathcal L} + I )f -\epsilon (f \cdot \nabla)f)\|.
$$
To estimate $\|{\mathcal L}^2 f \| $ it is necessary to have, at
least in the $y$ direction, control over six derivatives of $f_2$,
the second component of $f$.
 The reason is that
three derivatives of $f_2$ in the $y$ direction are necessary to
control ${\mathcal L} f$, due to the integral part of the operator
$\mathcal L$. Moreover, derivatives of different orders of the
perturbation may have different scales with respect to the
Reynolds number. These facts are shown in details in section
\ref{section5}. As already mentioned, clarifying and showing these
estimates in details is our aim, since it is a point that has been
overseen in previous works, and lead to mistakes about the
necessary assumptions on $f$.

We now continue the proof of stability. As mentioned above, since
$\| f \|^2_{H^6_m} = 1$, we can replace all the terms depending on
$f$ by an absolute constant and write inequality (\ref{ineq1}) as
\begin{eqnarray}
    \lefteqn{\displaystyle \int_0^{T}  \left( \| v \|^2_{\widetilde{H}} + \| v_t\|^2_{\widetilde{H}}\right) dt
         \leq
           C R^2
          + C R^2 \epsilon^2
           \int_0^{T} \| (v\cdot \nabla )v \|^2 dt } \hspace{3.3cm}\nonumber\\
    & & \displaystyle \mbox{} + C R^2 \epsilon^2
   \int_0^{T} \left( \| (v_t \cdot \nabla )v \|^2 + \| (v\cdot \nabla )v_t \|^2 \right) dt \label{eq22} \\
  & & \displaystyle \mbox{} +
 C R^2 \epsilon^2 \int_0^T \left( \| e^{-t}(v\cdot \nabla) f\|^2 + \| e^{-t}(f\cdot \nabla) v \|^2 \right) dt \nonumber\\
  & & \displaystyle \mbox{}
       + C R^2 \epsilon^2 \int_0^T \left( \| e^{-t}(v_t \cdot \nabla) f \|^2
        + \| e^{-t}(f\cdot \nabla) v_t \|^2 \right) dt.\nonumber
   \end{eqnarray}
From now on, we fix the constant $C$. To finish the proof, we use
the following Lemma, which is proved later:
\begin{lemma}\label{lema1}
   There exists $\epsilon_0>0$, $\epsilon_0 = {\mathcal O}(R^{-3})$,
   such that if $0 \leq \epsilon < \epsilon_0$ then
   \begin{equation}\label{eqlemma}
      \displaystyle \int_0^{T}  \left ( \| v(\cdot ,t) \|^2_{\widetilde{H}}
+ \| v_t(\cdot ,t) \|^2_{\widetilde{H}}\right ) dt
         <  2 C R^2 \; , \; \forall \: T \geq 0 .
   \end{equation}
\end{lemma}
Now, using (\ref{eq4}) and a simple one dimensional sobolev
inequality, we have
$$
    \displaystyle \max_{a \leq t \leq b} |v(\cdot , t) |_{\infty}^2 \leq
    \widetilde{C} R \max_{a \leq t \leq b} \|v(\cdot , t)
\|_{\widetilde{H}}^2 \leq
    \widetilde{C} R \left ( 1 + \frac{1}{b-a}\right )
\int_a^b \|v(\cdot ,t)\|_{\widetilde{H}}^2 dt
    + \int_a^b \|v_t (\cdot ,t)\|_{\widetilde{H}}^2 dt.
$$
This implies
\begin{equation} \label{eq37}
    \displaystyle \sup_{a \leq t } |v(\cdot , t) |_{\infty}^2 \leq
    \widetilde{C} R
    \int_a^{\infty} \left( \|v(\cdot ,t)\|_{\widetilde{H}}^2
    + \|v_t (\cdot ,t)\|_{\widetilde{H}}^2\right) dt .
\end{equation}
Note that in view of Lemma \ref{lema1}, the right hand side of
inequality (\ref{eq37}) is finite. Letting $a\rightarrow \infty$
in (\ref{eq37}), we have that $\displaystyle \lim_{t\rightarrow
\infty} |v(\cdot , t) |_{\infty}^2 = 0$, which proves the theorem.

\paragraph{Proof of Lemma \ref{lema1}:}
First, note that for $T > 0$ small enough, inequality
({\ref{eqlemma}) obviously holds. Now, suppose it does not hold
for all $T\geq 0$, that is, there exists $T_0 > 0$ such that
\begin{equation} \label{eq23}
     \displaystyle \int_0^{T_0}  \left ( \| v(\cdot ,t) \|^2_{\widetilde{H}}
+ \| v_t(\cdot ,t) \|^2_{\widetilde{H}}\right ) dt
         =  2 C R^2 .
\end{equation}
Using (\ref{eq22}), we have:
\begin{equation}
\begin{array}{l} \label{eq24}
     \displaystyle 2 C R^2 =
     \int_0^{T_0}  \left( \| v \|^2_{\widetilde{H}} + \| v_t \|^2_{\widetilde{H}}\right) dt
     \leq
     C R^2  \\
     \displaystyle \mbox{} + C R^2 \epsilon^2 \int_0^{T_0} \left( \|  (v\cdot \nabla )v \|^2 +
     \| (v_t \cdot \nabla )v \|^2 +
     \| (v\cdot \nabla )v_t \|^2 \right) dt \\
     \displaystyle \mbox{} +
     \int_0^{T_0} \left( \| e^{-t}(v\cdot \nabla) f\|^2 + \| e^{-t}(f\cdot \nabla) v \|^2
       + \| e^{-t}(v_t \cdot \nabla) f \|^2
        + \| e^{-t}(f\cdot \nabla) v_t \|^2 \right) dt .
\end{array}
\end{equation}
We now estimate the integrands on the right hand side of
inequality (\ref{eq24}) by the integral on its left hand side. To
this end, we will use the inequalities (\ref{eq4}) and
 \begin{equation} \label{eq26}
    \displaystyle \max_{0\leq t \leq T_0} \|v(\cdot , t)\|^2_{\widetilde{H}} \leq
          \int_0^{T_0}   \left ( \|v(\cdot , t)\|^2_{\widetilde{H}}
+ \|v_t (\cdot , t)\|^2_{\widetilde{H}}\right ) dt = 2CR^2.
 \end{equation}
Since $\displaystyle \| v \|_{\widetilde{H}}^2 = \|v\|^2 +
\frac{1}{R}\|Dv\|^2 + \frac{1}{R^2}\|v_{xy}\|^2$ , we have
$\|Dv\|^2 \leq R \|v\|^2_{\widetilde{H}}$. Therefore, for each $0
\leq t \leq T_0$:
\begin{eqnarray}
   \displaystyle \|  \left \{ (v\cdot \nabla )v\right \} (\cdot ,t)\|^2 & \leq &
   |v(\cdot , t)|_{\infty}^2 \|Dv(\cdot , t)\|^2 \label{eq27} \\
   \displaystyle & \leq &
   \left ( \widetilde{C} R \|v(\cdot ,t)\|^2_{\widetilde{H}} \right )
   \left ( R \|v(\cdot ,t)\|^2_{\widetilde{H}} \right )
   \leq 2\widetilde{C}C R^4 \|v( \cdot ,t)\|^2_{\widetilde{H}}, \nonumber
\end{eqnarray}
\begin{eqnarray}
   \displaystyle \|  \left \{ (v_t \cdot \nabla )v\right \} (\cdot ,t)\|^2 & \leq &
   |v_t (\cdot , t)|_{\infty}^2 \|Dv(\cdot , t)\|^2 \label{eq28} \\
   & \leq &
   \left ( \widetilde{C} R \|v_t (\cdot ,t)\|^2_{\widetilde{H}} \right )
   \left ( R \|v(\cdot ,t)\|^2_{\widetilde{H}} \right ) \leq
   2\widetilde{C} C R^{4} \|v_t ( \cdot ,t)\|^2_{\widetilde{H}} , \nonumber
\end{eqnarray}
\begin{eqnarray}
   \displaystyle \|  \left \{ (v\cdot \nabla )v_t\right \} (\cdot ,t)\|^2 & \leq &
   |v(\cdot , t)|_{\infty}^2 \|Dv_t (\cdot , t)\|^2 \label{eq29} \\
   & \leq & \left ( \widetilde{C} R \|v(\cdot ,t)\|^2_{\widetilde{H}} \right )
   \left ( R \|v_t (\cdot ,t)\|^2_{\widetilde{H}} \right ) \leq
   2\widetilde{C} C R^{4} \|v_t ( \cdot ,t)\|^2_{\widetilde{H}}, \nonumber
\end{eqnarray}
\begin{eqnarray} \label{eq30}
   \displaystyle \|  e^{-t} \left\{ (v\cdot \nabla )f \right \}(\cdot ,t)\|^2 \leq
   e^{-2t}|v(\cdot , t)|_{\infty}^2 \| Df \|^2 \leq
   \widetilde{C} R \|v(\cdot ,t)\|^2_{\widetilde{H}},
\end{eqnarray}
\begin{eqnarray} \label{eq31}
   \displaystyle \|  e^{-t} \left \{ (f\cdot \nabla )v \right \}(\cdot ,t)\|^2 \leq
   e^{-2t}|f|_{\infty}^2 \| Dv(\cdot ,t) \|^2 \leq
    R \|v(\cdot ,t)\|^2_{\widetilde{H}},
\end{eqnarray}
\begin{eqnarray} \label{eq32}
   \displaystyle \|  e^{-t} \left \{ (v_t \cdot \nabla )f \right \} (\cdot ,t)\|^2 \leq
   e^{-2t}|v_t (\cdot , t)|_{\infty}^2 \| Df \|^2 \leq
   \widetilde{C} R \|v_t (\cdot ,t)\|^2_{\widetilde{H}},
\end{eqnarray}
\begin{eqnarray} \label{eq33}
   \displaystyle \|  e^{-t} \left \{ (f\cdot \nabla )v_t \right \} (\cdot ,t)\|^2 \leq
   e^{-2t}|f|_{\infty}^2 \| Dv_t (\cdot ,t) \|^2 \leq
   R \|v_t (\cdot ,t)\|^2_{\widetilde{H}} .
\end{eqnarray}
Applying (\ref{eq27}), (\ref{eq28}), (\ref{eq29}), (\ref{eq30}),
(\ref{eq31}), (\ref{eq32}), (\ref{eq33}) to (\ref{eq24}) gives
\begin{equation} \label{eq34}
   \begin{array}{l}
     \displaystyle 2 C R^2
     \leq
     C R^2
     +  C R^2 \epsilon^2 \left \{
      6 \widetilde{C} C R^4 \int_0^{T_0} \left( \|  v(\cdot ,t) \|^2_{\tilde{H}} +
     \| v_t (\cdot ,t) \|^2_{\tilde{H}} \right) dt \right\}  \\
     = C R^2 + C R^2 \epsilon^2 \left \{
      12 \widetilde{C} C^2 R^6 \right \}.
   \end{array}
\end{equation}
This implies
\begin{equation} \label{eq35}
   \displaystyle 1 \leq 12 \widetilde{C}
        C^2 R^6 \epsilon^2 ,
\end{equation}
which is equivalent to
\begin{equation}
 \epsilon \geq
        \frac{1}{CR^3\sqrt{12\widetilde{C}}} =
         \frac{1}{K R^3}
\end{equation}
where $\displaystyle K := C\sqrt{12\widetilde{C}}$. Therefore, if
$\displaystyle \epsilon < \frac{1}{K R^3}$, equality (\ref{eq23})
never holds. This proves the Lemma.

\section{Estimates for the pressure terms} \label{section5}

We now turn our attention to our main objective, which is to show
why one needs to assume control over at least six derivatives of
$f_2$ in the $y$ direction to apply the method above. This
necessity follows from estimates (\ref{eqp1}) and (\ref{eqp2}) for
the pressure terms $p_1(x,y,t)$ and $p_2(x,y,t)$. So, we begin by
showing these estimates.

\begin{theorem}\label{theorempressure}
 If $p_1(x,y,t)$, $p_2(x,y,t)$ are the solutions of
\begin{equation} \label{eqb1}
    \left \{
\begin{array}{l}
        \displaystyle \Delta p_1 = -\nabla \cdot ((u\cdot
\nabla )U) -  \nabla \cdot ((U\cdot \nabla )u) \\
        \displaystyle
p_{1y}(x,0,t) = \frac{1}{R} u_{2yy}(x,0,t) \vspace{.1cm}\\
\displaystyle p_{1y}(x,1,t) = \frac{1}{R} u_{2yy}(x,1,t)
     \end{array}
\right.
\end{equation}
and
\begin{equation}\label{eqb2}
\left \{ \begin{array}{l}
        \displaystyle \Delta p_2 = -\nabla \cdot
((u\cdot \nabla )u) \\
        \displaystyle p_{2y}(x,0,t) =
0\vspace{.1cm}\\
        \displaystyle p_{2y}(x,1,t) = 0 ,
     \end{array}
\right.
\end{equation}
then
\begin{eqnarray}
\|\nabla p_1 (\cdot , \cdot , t )\|^2 & \leq & C \left( \|u (\cdot
, \cdot , t )\|_{H^1}^2 + \frac{1}{R^2} \| u_{2yy}(\cdot , \cdot ,
t ) \| + \frac{1}{R^2} \| u_{2yyy} (\cdot , \cdot , t )\| \right)
    \label{eqb3}\\
\|\nabla p_2 (\cdot , \cdot , t )\|^2 & \leq & \|( u\cdot \nabla)
u(\cdot , \cdot , t )\|^2 . \label{eqb4}
\end{eqnarray}
for all $t\geq 0$, where $C$ is an absolute constant.
\end{theorem}
Note that the inequalities above are for norms with respect to the
space variables $x$ and $y$. Therefore, to simplify the notation,
we prove them for functions depending only on these variables. As
before, $\Omega = [0,1]\times [0,1]$, and all norms are over
$\Omega$. For clarity of the presentation, we separate the proof
of Theorem \ref{theorempressure} into two Lemmas.
\begin{lemma} \label{lemab1}
Let $g : \mathbb R \times [0 , 1] \rightarrow \mathbb R^2$,
$g(x,y) = (g_1(x,y) , g_2(x,y))$, be a $C^\infty$ function
satisfying
\begin{eqnarray}
& & g(x,1) = g( x, 0) = (0,0) \label{eqb5}\\
& & g(x,y) = g(x+1,y) \, \forall \, x \in \mathbb R.\label{eqb6}
\end{eqnarray}
If $h: \mathbb R\times [0 , 1] \rightarrow \mathbb R$ is the
solution of
\begin{equation}\label{eqb7}
\left \{ \begin{array}{l}
        \displaystyle \Delta h = \nabla \cdot g \\
        \displaystyle h_{y}(x,0) = 0 \vspace{.1cm}\\
        \displaystyle h_{y}(x,1) = 0 \\
        h(x,y) = h(x+1,y) ,
     \end{array}
\right.
\end{equation}
then
\begin{equation}\label{eqb8}
\|\nabla h \|^2  \leq  \|g\|^2 .
\end{equation}
\end{lemma}
\paragraph{Proof:} If $\Delta h = \nabla \cdot g$, then
$$
 \displaystyle \int_\Omega (h_{xx} + h_{yy})h \, dxdy = \int_\Omega (g_{1x} + g_{2y}) h \, dxdy.
$$
Through integration by parts,
\begin{eqnarray*}
  \lefteqn{- \displaystyle \int_\Omega (h_x^2 + h_y^2) \, dx dy +  \int_{\partial\Omega} h_x h \nu^x \, dS
  +  \int_{\partial\Omega} h_y h \nu^y \, dS  = } \\
  & & \mbox{}  - \int_{\Omega} g_1 h_x\, dx dy - \int_{\Omega} g_2 h_y \, dx dy
  +  \int_{\partial\Omega} g_1 h \nu^x \, dS + \int_{\partial\Omega} g_2 h \nu^y \, dS ,
\end{eqnarray*}
where $\nu^x$ and $\nu^y$ denote the components of the outer
normal to $\partial \Omega$ in the $x$
 and $y$ directions respectively.
From the conditions satisfied by $h$ and $g$ at the boundary, the
boundary integrals above vanish. Then,
$$
 \|\nabla h\|^2 = \displaystyle \int_\Omega (h_x^2 + h_y^2) \, dx dy =
 \int_{\Omega} g_1 h_x\, dx dy + \int_{\Omega} g_2 h_y \, dx dy .
$$
Using the Cauchy-Schwarz inequality,
\begin{eqnarray*}
\displaystyle  \|\nabla h\|^2 \leq \|g_1\| \|h_x\| + \|g_2\|
\|h_y\| & \leq &
 \frac{1}{2}\|g_1\|^2 + \frac{1}{2}\|h_x\|^2 + \frac{1}{2}\|g_2\|^2 + \frac{1}{2}\|h_y\|^2 \\
 & = &
 \frac{1}{2} \| g \|^2 + \frac{1}{2} \|\nabla h\|^2 .
\end{eqnarray*}
This implies the desired estimate
$$
  \|\nabla h \|^2 \leq \| g \|^2.
$$

The Lemma above gives inequality (\ref{eqb4}) for $p_2(x,y,t)$,
the solution of (\ref{eqb2}).

As mentioned before, the estimates to be proved do not depend on
the variable $t$. Therefore, we write simply $u(x,y) = (u_1 (x,y)
, u_2(x,y))$, for $u$ the solution of problem (\ref{eq5}). We
remind the reader that $U(x,y) = (y,0)$. We prove the following
lemma, completing the proof of Theorem \ref{theorempressure}:

\begin{lemma} \label{lemab2}
If $h: \mathbb R\times [0 , 1] \rightarrow \mathbb R$ is the
solution of
\begin{equation}\label{eqb9}
\left \{ \begin{array}{l}
        \displaystyle \Delta h = - \nabla \cdot (u \cdot \nabla )U - \nabla \cdot (U \cdot \nabla ) u \\
        \displaystyle h_{y}(x,0) =  \frac{1}{R} u_{2yy}(x,0) \vspace{.1cm}\\
        \displaystyle h_{y}(x,1) =  \frac{1}{R} u_{2yy}(x,1) \vspace{.1cm}\\
        h(x,y) = h(x+1,y) ,
     \end{array}
\right.
\end{equation}
then
\begin{equation}\label{eqb10}
\|\nabla h \|^2  \leq C \left( \|u\|_{H^1}^2 +
\frac{1}{R^2}\|u_{2yy}\|^2 + \frac{1}{R^2}\| u_{2yyy} \|^2\right).
\end{equation}
\end{lemma}
\paragraph{Proof:} We begin by noting that if
$h_1$, $h_2$ are the solutions of
\begin{equation}\label{eqb11}
\left \{ \begin{array}{l}
        \displaystyle \Delta h_1 = - \nabla \cdot (u \cdot \nabla )U - \nabla \cdot (U \cdot \nabla ) u \\
        \displaystyle h_{1y}(x,0) =  0\\
        \displaystyle h_{1y}(x,1) =  0\\
        h_1(x,y) = h_1(x+1,y) ,
     \end{array}
\right.
\end{equation}
and
\begin{equation}\label{eqb12}
\left \{ \begin{array}{l}
        \displaystyle \Delta h_2 = 0\\
        \displaystyle h_{2y}(x,0) =  \frac{1}{R} u_{2yy}(x,0) \vspace{.1cm}\\
        \displaystyle h_{2y}(x,1) =  \frac{1}{R} u_{2yy}(x,1) \vspace{.1cm}\\
        h_2(x,y) = h_2(x+1,y) ,
     \end{array}
\right.
\end{equation}
then $h = h_1 + h_2 $ is the solution of problem (\ref{eqb9}).
Therefore, to prove (\ref{eqb10}), it is sufficient to prove
estimates for $h_1$ and $h_2$, solutions of (\ref{eqb11}) and
(\ref{eqb12}) respectively. For $h_1$, lemma \ref{lemab1} implies
\begin{equation}\label{eqb13}
  \| \nabla h_1 \|^2 \leq  \|(u \cdot \nabla )U + (U \cdot \nabla ) u\|^2 \leq 2 \|u\|^2 + 2 \| D u\|^2 =
  2\|u \|^2_{H^1}.
\end{equation}
To prove the estimates for $h_2$, expand in a Fourier series in
the $x$ direction. The Fourier coefficients $\hh_2 (k , y)$
satisfy
\begin{eqnarray}
        \displaystyle -k^2 \hh_2 + \widehat{h}_2^{\prime\prime} & = & 0 \label{eqb14}\\
        \displaystyle \hh_{2y}(k,0) & = & \frac{1}{R} \hu_{2yy}(k,0) \label{eqb15}\\
        \displaystyle \hh_{2y}(k,1) & = &  \frac{1}{R} \hu_{2yy}(k,1) , \label{eqb16}
\end{eqnarray}
where $^\prime$ denotes the derivative with respect to $y$.
Consider first the case $k\neq 0$. To simplify notation, let
\begin{eqnarray*}
\alpha_k & := & \frac{1}{R} \hu_{2yy}(k,0) \vspace{.1cm} \\
\beta_k & := & \frac{1}{R} \hu_{2yy}(k,1).
\end{eqnarray*}
Using a one-dimensional sobolev type inequality, we estimate
$\alpha_k$ and $\beta_k$ by
\begin{eqnarray}
 & & \displaystyle |\alpha_k |^2 \leq
     \displaystyle \frac{1}{R^2}\max_{0\leq y\leq1} | \hu_{2yy} (k , y) |^2
     \leq \frac{C}{R^2}\left( \|\hu_{2yy} (k , \cdot)\|^2 +
     \| \hu_{2yyy} (k , \cdot)\|^2\right) \label{eqb163}\\
 & & \displaystyle |\beta_k |^2 \leq
     \displaystyle \frac{1}{R^2}\max_{0\leq y\leq1} | \hu_{2yy} (k , y) |^2
     \leq \frac{C}{R^2}\left( \|\hu_{2yy} (k , \cdot)\|^2 +
     \| \hu_{2yyy} (k , \cdot)\|^2\right),\label{eqb165}
\end{eqnarray}
where $C$ is an absolute constant. As before, we keep the notation
simple by using $C$ to represent any absolute constant, whose
value can possibly change for different inequalities.\\
The general solution of the differential equation (\ref{eqb14}) is
\begin{equation}\label{eqb17}
  \widehat{h}_2 (k , y) = a_k e^{|k|(y-1)} + b_k e^{-|k|y} .
\end{equation}
Imposing the boundary conditions (\ref{eqb15}) and (\ref{eqb16}),
we determine the coefficients $a_k$ and $b_k$:
$$
 \begin{array}{lcr}
  \displaystyle a_k = \frac{e^{|k|}}{e^{|k|} - e^{-|k|}} \frac{\beta_k}{|k|} - \frac{1}{e^{|k|} - e^{-|k|}} \frac{\alpha_k}{|k|} & ; &
  \displaystyle b_k = \frac{1}{e^{|k|} - e^{-|k|}} \frac{\beta_k}{|k|} - \frac{e^{|k|}}{e^{|k|} - e^{-|k|}} \frac{\alpha_k}{|k|} .
 \end{array}
$$
Therefore,
\begin{eqnarray}
 \displaystyle |a_k|^2 & \leq & C \left( \frac{|\alpha^2_k|}{k^2} + \frac{|\beta_k|^2}{k^2}\right) \label{eqb18}\\
 \displaystyle |b_k|^2 & \leq & C \left( \frac{|\alpha^2_k|}{k^2} + \frac{|\beta_k|^2}{k^2}\right). \label{eqb19}
 \end{eqnarray}
Using (\ref{eqb163}), (\ref{eqb165}), (\ref{eqb17}),
({\ref{eqb18}), (\ref{eqb19}), we have
 \begin{eqnarray}
   \displaystyle k^2 \| \widehat{h}_2 (k , \cdot)\|^2 & \leq & C k^2 ( |a_k|^2 \| e^{|k|(y-1)}\|^2 + |b_k|^2 \| e^{-|k|y}\|^2 )
     \leq  \displaystyle \frac{C}{|k|} (|\alpha_k |^2 + |\beta_k |^2 ) \nonumber \\
    & \leq & C (|\alpha_k |^2 + |\beta_k |^2 ) \leq
    \displaystyle \frac{C}{R^2}\left( \|\hu_{2yy} (k , \cdot)\|^2 +
     \| \hu_{2yyy} (k , \cdot)\|^2\right) \label{eqb20}
 \end{eqnarray}
and
 \begin{eqnarray}
   \displaystyle \| \widehat{h}_{2y} (k , \cdot)\|^2 & \leq & C k^2 ( |a_k|^2 \| e^{|k|(y-1)}\|^2 + |b_k|^2 \| e^{-|k|y}\|^2 )
     \leq  \displaystyle \frac{C}{|k|} (|\alpha_k |^2 + |\beta_k |^2 ) \nonumber\\
    & \leq & C (|\alpha_k |^2 + |\beta_k |^2 ) \leq \displaystyle \frac{C}{R^2}\left( \|\hu_{2yy} (k , \cdot)\|^2 +
     \| \hu_{2yyy} (k , \cdot)\|^2\right). \label{eqb21}
 \end{eqnarray}
For the $k = 0$ mode, we solve equation (\ref{eqb14}) under
boundary conditions (\ref{eqb15}), (\ref{eqb16}) directly. Note
that the divergence free condition satisfied by $u$ assures that
this problem is solvable, as expected.
 Differentiating the solution with respect to $y$, we get
$$
\widehat{h}_{2y} (0 , y) = \displaystyle \frac{1}{R} \hu_{2yy} (0
, 1).
$$
Then
\begin{equation}\label{eqb22}
 \| \widehat{h}_{2y} (0 , \cdot)\|^2 \leq
 \frac{1}{R^2}\max_{0\leq y\leq1} | \hu_{2yy} (0 , y) |^2
\leq \frac{C}{R^2}\left( \|\hu_{2yy} (0 , \cdot)\|^2 +
     \| \hu_{2yyy} (0 , \cdot)\|^2\right).
\end{equation}
Using (\ref{eqb20}), (\ref{eqb21}) and (\ref{eqb22}),
\begin{equation}\label{eqb23}
\begin{array}{crl}
\displaystyle \|\nabla h_2 \|^2 & = & \displaystyle\sum_{k\in
\mathbb Z}
\left( k^2 \displaystyle\|\widehat{h}_2 (k , \cdot) \|^2 + \|\widehat{h}_{2y} ( k , \cdot) \|^2\right)\vspace{.1cm}\\
 & \leq & \displaystyle \frac{C}{R^2}\sum_{k\in \mathbb Z}\left( \|\hu_{2yy} (k , \cdot)\|^2 +
     \| \hu_{2yyy} (k , \cdot)\|^2\right) \vspace{.1cm}\\
 & = & \displaystyle \frac{C}{R^2} \left( \|u_{2yy}\|^2 + \| u_{2yyy} \|^2\right).
\end{array}
\end{equation}
Therefore, by (\ref{eqb13}) and (\ref{eqb23}), we have that $h
=h_1 + h_2$, solution of (\ref{eqb9}), satisfies
\begin{equation}\label{eqb24}
  \|\nabla h \|^2 \leq C \left( \|u\|_{H^1}^2 + \frac{1}{R^2}\|u_{2yy}\|^2 + \frac{1}{R^2}\| u_{2yyy} \|^2\right)
  .
\end{equation}
This finishes the proof of the Lemma. This Lemma applied to
$p_1(x,y,t)$ completes the proof of theorem \ref{theorempressure}.

In section \ref{section4}, we used that the terms depending only
on $f$ on the right hand side of (\ref{ineq1}) can be estimated by
$C \| f \|^2_{H^6_m}$. This fact is a consequence of Theorem
\ref{theorempressure}. In fact, it is easy to see that the term
that requires highest derivatives of $f$ to be estimated is $\|
{\mathcal L}^2 f\|^2$. Therefore, we should bound this term in a
sharp way, that is, estimate it using a norm for $f$ that involves
derivatives of lowest possible order. We show this norm to be $\|
f \|_{H^6_m}$.

Recall that
 \begin{equation}
\displaystyle {\mathcal L} f = \frac{1}{R}\Delta f - (f\cdot
\nabla )U - (U\cdot \nabla )f - \nabla p ,
\end{equation}
where $U = (y , 0)$, and $p$ is the solution of
\begin{equation}
    \left \{ \begin{array}{l}
        \displaystyle \Delta p = -\nabla \cdot ((f\cdot \nabla )U) -  \nabla \cdot ((U\cdot \nabla )f) \\
        \displaystyle p_{y}(x,0) = \frac{1}{R} f_{2yy}(x,0) \vspace{.1cm}\\
        \displaystyle p_{y}(x,1) = \frac{1}{R} f_{2yy}(x,1) .
     \end{array} \right.
\end{equation}
Using Theorem \ref{theorempressure}, we get that
\begin{equation} \label{estimate1}
 \| {\mathcal L} f \|^2 \leq C \left( \| f \|^2_{H^1} + \frac{1}{R^2} \| f_{xx} \|^2
 + \frac{1}{R^2} \| f_{yy} \|^2
 + \frac{1}{R^2} \| f_{2yyy} \|^2 \right).
\end{equation}
Therefore,
\begin{equation} \label{estimate2}
 \| {\mathcal L}^2 f \|^2 \leq C \left( \| {\mathcal L}f \|^2_{H^1} +
  \frac{1}{R^2} \|{\mathcal L} f_{xx} \|^2 + \frac{1}{R^2} \|{\mathcal L} f_{yy} \|^2
 + \frac{1}{R^2} \| ({\mathcal L}f)_{2yyy} \|^2 \right).
\end{equation}
Straightforward computations show that
\begin{eqnarray*}
 \| {\mathcal L}f \|^2_{H^1} & \leq & C \left( \| f \|^2_{H^2}
 + \frac{1}{R^2} \| D^3 f\|^2 + \frac{1}{R^2} \| D^4 f \|^2 \right) \\
 \frac{1}{R^2} \|{\mathcal L} f_{xx} \|^2 & \leq & \frac{C}{R^2} \left( \| f \|^2_{H^2} +
 \| D^3 f \|^2 + \frac{1}{R^2} \| D^4 f \|^2 + \frac{1}{R^2} \|f_{2yyyxx} \|^2 \right) \\
 \frac{1}{R^2} \|{\mathcal L} f_{yy} \|^2 & \leq & \frac{C}{R^2} \left( \| f \|^2_{H^2} +
 \| D^3 f \|^2 +
 \frac{1}{R^2} \| D^4 f \|^2 + \frac{1}{R^2} \|f_{2yyyyy} \|^2 \right) \\
\frac{1}{R^2} \| ({\mathcal L}f)_{2yyy} \|^2 & \leq &
\frac{C}{R^2} \left(
 \| f_{2yyy} \|^2 + \| f_{2xyyy}\|^2 \right) \\ & &  \mbox{}+ \frac{C}{R^4} \left( \| f_{2yyyxx} \|^2 +
 \| f_{2yyyyy} \|^2 +  \|f_{2yyyyyy} \|^2 \right) .
\end{eqnarray*}
The inequalities above, together with (\ref{estimate2}), imply
\begin{eqnarray}
 \| {\mathcal L}^2 f \|^2 & \leq & C \left( \| f \|^2_{H^2}
 + \frac{1}{R^2} \| D^3 f\|^2 + \frac{1}{R^2} \| D^4 f \|^2 +
\frac{1}{R^4}\| f_{2xxyyy}\|^2 \right.\nonumber \\
& & \mbox{} + \left. \frac{1}{R^4}\| f_{2yyyyy}\|^2 +
\frac{1}{R^4}\| f_{2yyyyyy}\|^2 \right) = C \| f \|^2_{H^6_m} .
\end{eqnarray}
Note that this estimate is sharp in view of  theorem
\ref{theorempressure} and the arguments above: one needs at least
six derivatives of $f_2$ to be able to estimate $\| {\mathcal L}^2
f \|^2$.

\section{conclusions}

As showed in section \ref{section4}, perturbations with norm $\|
\cdot \|_{H^6_m}$ of order $R^{-3}$ decay with time. Note that
this also shows that even though one needs control over
derivatives of high order of the perturbation, there are different
scales for those derivatives. In fact, if
\begin{eqnarray}
 \displaystyle \|f\|^2_{H^6_m} & = & \|f\|^2_{H^2} + \frac{1}{R^2}\|D^3 f\|^2 +
 \frac{1}{R^2}\|D^4 f\|^2 +
 \frac{1}{R^4}\| f_{2xxyyy}\|^2 \nonumber \\& & \mbox{} + \frac{1}{R^4}\| f_{2yyyyy}\|^2 +
\frac{1}{R^4}\| f_{2yyyyyy}\|^2 = {\mathcal O}(R^{-6}),
\end{eqnarray}
then $\| f \|_{H^2} $ is of order $R^{-3}$, $\|D^3 f \|$ and
$\|D^4 f \|$ are of order $R^{-2}$, and $\| f_{2xxyyy} \|$,
$\|f_{2yyyyy} \|$ and $\|f_{2yyyyyy} \|$ are of order $R^{-1}$.

In applications, control over this many derivatives is too
restrictive. One possible way to avoid this requirement is to
incorporate a smoothing property of the system to the argument. We
hope to address this question in the future.


\begin{thebibliography}{99}

\bibitem{BS} P. Braz e Silva,
Resolvent estimates for 2 dimensional perturbations of plane
Couette flow, Electron. J. Diff. Eqns., (2002) 92 (2002), 1--15.

\bibitem{KL} H.-O. Kreiss and J. Lorenz, "Initial-Boundary Value Problems
and the Navier-Stokes Equations", Pure and Applied Mathematics 136
Academic Press, 1989.

\bibitem{KL2} H.-O.Kreiss and J. Lorenz, Stability for time dependent differential
equations, Acta Numer. (7)(1998), 203-285, Cambridge University
Press, Cambridge.

\bibitem{KL3} H.-O. Kreiss and J. Lorenz, Resolvent estimates and quantification of
nonlinear stability, Acta Math. Sin. (Engl. Ser.) (16) 1 (2000),
1-20

\bibitem{K} G. Kreiss, A. Lundbladh and D.S. Henningson,
Bounds for threshold amplitudes in subcritical shear flows, J.
Fluid Mech., (270)(1994), 175--198.

\bibitem{L1} M. Liefvendahl and G. Kreiss, Bounds for the
threshold amplitude for plane Couette flow, J. Nonlinear Math.
Phys. (9) 3 (2002), 311--324.

\bibitem{L2} M. Liefvendahl and G. Kreiss, Analytical and
numerical investigation of the resolvent for plane Couette flow,
SIAM J. Appl. Math. (63) 3 (2003), 801--817.

\bibitem{RE} S.C. Reddy and D.S. Henningson, Energy growth in
viscous channel flows, J. Fluid Mech. (252) (1993), 209--238.

\bibitem{R} V.A. Romanov, Stability of plane-parallel Couette
flow, Functional Anal. Applics. (7) (1973), 137--146.

\bibitem{T} L.N. Trefethen, A.E. Trefethen, S.C. Reddy and T.A. Driscoll,
Hydrodynamic Stability Without Eigenvalues, Science (261) (1993),
578-584.


\end{thebibliography}
\end{document}